\documentclass[11pt,twoside]{amsart}

\usepackage{amssymb, palatino,graphicx, hyperref}
\usepackage[normalem]{ulem}
\usepackage[all]{xy}
\usepackage{xcolor}
\usepackage{ifthen}

\date{\today}

\newtheorem{dfn}{Definition}[section]
\newtheorem{thm}[dfn]{Theorem}
\newtheorem{crl}[dfn]{Corollary}
\newtheorem{prop}[dfn]{Proposition}
\newtheorem{rem}[dfn]{Remark}

\def\P{\ensuremath{\mathbb P}}
\def\N{\ensuremath{\mathbb N}}
\def\Q{\ensuremath{\mathbb Q}}
\def\R{\ensuremath{\mathbb R}}
\def\Z{\ensuremath{\mathbb Z}}

\def\C{\ensuremath{\mathbb C}}

\def\t{\mathfrak{t}}

\def\CP1{\ensuremath{\mathbb C \mathbb P^1}}

\def\Cn{\ensuremath{\C^n}}

\def\Pn-1{\ensuremath{\P^{n-1}}}

\newcommand{\h}{\ensuremath{ \mathfrak{h} }}

\def\cf{{\em cf. }}

\newcommand{\setofst}[2]
{ \ensuremath{ \left\{\, #1\ \left\vert\ #2 \right.\right\}
            }
}

\newcommand{\horiz}{\ensuremath $\rotatebox{90}{$\vert$}$}

\definecolor{red}{rgb}{.6,0,0}
\definecolor{green}{rgb}{0,.6,0}
\definecolor{darkgreen}{rgb}{0,0.3,0}
\definecolor{purple}{rgb}{0.5,0,0.5}
\definecolor{darkblue}{rgb}{0,0,0.7}
\definecolor{greenblue}{rgb}{0,0.4,0.5}
\definecolor{myblue}{HTML}{1685d6}
\definecolor{mypurple}{HTML}{b82e6c}

\newboolean{draft}

\newcommand{\cmt}[1]
{\ifthenelse {\boolean{draft}}
{{\sc \tiny \color{red} #1}}
{}}

\newcommand{\newbb}[1]
{\ifthenelse {\boolean{draft}}
{{\color{darkblue} #1}}
{#1}}

\newcommand{\newbbb}[1]
{\ifthenelse {\boolean{draft}}
{{\color{greenblue} #1}}
{#1}}

\newcommand{\nopost}[1]
{\ifthenelse {\boolean{draft}}
{{\color{cyan} #1}}
{}}

\newcommand{\maynopost}[1]
{\ifthenelse {\boolean{draft}}
{{\color{purple} #1}}
{}}

\newcommand{\margincmt}[1]
{\ifthenelse {\boolean{draft}}
{\marginpar{{\sc \tiny \color{red} #1}}}
{}}
\newcommand{\inred}[1]
{\ifthenelse{\boolean{draft}}{{\color{red} #1}}{#1}}

\newcommand{\new}[1]
{\ifthenelse {\boolean{draft}}
{{\color{green} #1}}
{#1}}

\newcommand{\neww}[1]
{\ifthenelse {\boolean{draft}}
{{\color{darkgreen} #1}}
{#1}}

\newcommand{\newb}[1]
{\ifthenelse {\boolean{draft}}
{{\color{blue} #1}}
{#1}}

\newcommand{\del}[1]
{\ifthenelse {\boolean{draft}}
{{\color{magenta} #1}}
{}}

\newboolean{details_on}

\newcommand{\details}[1]
{\ifthenelse {\boolean{details_on}}
{{\color{darkgreen} \tiny #1}}
{}}

\title{Simplicial toric varieties as leaf spaces}
\author{Fiammetta Battaglia and Dan Zaffran}
\thanks{}

\newcommand\vz{\ensuremath{\underline{z}}}

\newcommand{\bT}{\ensuremath{\mathcal{T} } } 

\newcommand{\calF}{\ensuremath{ \mathcal{F} } }

\newcommand{\lr}[1]{\ensuremath{ \Lambda^\R_{#1} } }

\DeclareMathOperator{\Rel}{Rel}

\setboolean{draft}{true}
\setboolean{details_on}{true}

\begin{document}
\maketitle

\begin{abstract} We present a summary of some results from our article \cite{bz1} and other recent results on the so-called LVMB manifolds. We emphasize some features by taking a different point of view. We present a simple variant of the Delzant construction, in which the group that is used to perform the symplectic reduction can be chosen of arbitrarily high dimension, and is always connected. 
\end{abstract}
\section*{Introduction}
At the workshop held in Rome in 16-20 November 2015, for the 60th birthday of Simon Salamon, the first author presented a joint article with Dan Zaffran \cite{bz1}. In this note we make a summary of various results, including previous and subsequent literature, concerning the relation between a special class of compact foliated complex manifolds, called LVMB manifolds, toric geometry, and convex geometry; we will simultaneously treat the rational and nonrational cases.
In particular we will dwelve on some aspects of the article \cite{bz1} that have been left aside in the published version  or that were not developed at the time. Among the last, a variant of the Delzant procedure. 

We are interested in the relationship between three different classes of objects; we first describe each of them briefly.
  
LVM manifolds originated in the context of dynamical systems. They form a large class of non-K\"ahler, compact, complex manifolds, introduced between 1997 and 2001 in works by Lopez de Medrano, Verjovsky, and Meersseman \cite{ldm,M}. Their construction was generalized by Bosio \cite{Bos}, who obtained a larger class, which we will refer to as LVMB manifolds.
The classical starting datum for the construction of an LVMB manifold is a set of holomorphic vector fields, inducing a
$\C^m$-action on $\C^n$, together with a choice of a saturated open subset of $\C^n$ for which the space of orbits
is a compact
complex manifold, interestingly non K\"ahler. Each LVMB manifold $N$ is also endowed with a holomorphic foliation 
${\calF}$, as shown in [16,17] for LVM manifolds and in [8] for LVMB manifolds.

The theory of toric varieties is by now classical. Simplicial toric varieties are algebraic manifolds with at most finite quotient singularities. They can be seen as compactifications of a torus $(\C^*)^n$. There are several reference texts for
toric varieties theory, among them \cite{Ful,coxtoric,oda}. 
The standard starting datum for the construction of a compact simplicial toric variety is a complete simplicial rational fan. 

A fan is a set of convex polyhedral cones having certain properties. Recall that a cone in a vector space is the set of nonnegative linear combinations of a finite number of vectors, that generate the cone. Each vector generates a nonnegative half line, called a ray of the cone. A cone is simplicial if it admits a set of linearly independent generators. Let $L$ be a lattice in the vector space $L\otimes_{\Z}\R$. A cone in $L\otimes_{\Z}\R$ is rational if each of its rays has nonempty intersection with $L$. A fan is simplicial if each of its cones is simplicial, and it is rational in $L\otimes_{\Z}\R$ if each of its cones is rational. In this work we will consider fans and other related convex objects.

Now, is there a relation between LVMB manifolds and toric varieties? Is there a relation between certain linear $\C^m$-actions on $\C^n$ and fans? What happens when the fan is nonrational? What is a nonrational fan? Is a fan the appropriate convex object? What are the similarities and differences between LVMB manifolds and toric manifolds? In the present article we try to give an answer to these questions.

Recall that a fan in a  vector space is complete if the union of its cones is the whole space; it is polytopal if it is the fan normal to a polytope. In \cite{MV} Meersseman and Verjovsky establish a precise relationship between LVM manifolds and compact simplicial {\em projective} toric varieties---that is, varieties associated to complete, simplicial, {\em polytopal} fans. They prove that the leaf space of an LVM manifold whose starting datum satisfies a further rationality condition---condition (K)---is a simplicial projective toric variety. Conversely, any simplicial projective toric variety can be obtained as leaf space of an LVM manifolds of that kind. In order to prove this last result they use Gale duality combined with symplectic reduction. 
 
But what happens when the fan is not polytopal or nonrational?

Cupit-Foutou and Zaffran in \cite{CZ} establish that the class of LVM manifolds is strictly included in the LVMB family. Battisti further proves in \cite{B} that an LVMB manifold is LVM if and only if the corresponding fan is polytopal. 
  
How can we deal with the nonrational case? In classical toric geometry, when one considers a rational fan in $\R^d$, there are two data that are usually taken for granted: a lattice, and, in the lattice, a set of primitive vectors, each of which a generator of a fan ray. In order to extend this setting to the nonrational case one needs to reconsider these data.
Let us illustrate how from our view-point. We introduce a convex object that allows to encode all of these data: a triangulated vector configuration. This is a pair $\{V,\bT\}$, where 
$V=(v_1,\ldots,v_n)$ is a configuration (ordered and allowing repetitions) of vectors in $\R^d$, and $\bT$ is a collection of subsets of $\{1,\ldots,n\}$ with suitable properties, called a triangulation of $V$. Consider for example a rational fan in a lattice $L$ with $h$ rays, and, for each ray, its primitive generator. Then, a corresponding triangulated vector configuration is a (non-unique) pair $\{V,\bT\}$ such that: $\text{Span}_\Z(V)=L$; the first $h$ vectors in $V$ are the selected generators of the $h$ fan rays; the triangulation $\bT$ carries the combinatorial information that determines the subcollections of $\{v_1,\ldots,v_n\}$ that generate all of the fan cones. Notice that $n>h$ may be needed, for example in case the set of primitive ray generators is not a generating set of $L$.
As an example, consider the rational simplicial fan in $\R^2=\Z^2\otimes_{\Z}\R$ drawn in the picture. Its associated
toric variety is $\C P^1\times\C P^1$ blown up at one point. A corresponding triangulated vector configuration is $\{V,\bT\}$, with $V=(\ (1,0),(0,1),(-1,0),(-1,1),(0,-1)\ )$
and $\bT$ the triangulation whose maximal simplices are $\{\{1,2\},\{2,3\},\{3,4\},\{4,5\},\{5,1\}\}$; here $h=n=5$.  
\begin{center}
\def\svgwidth{0.5\textwidth}
\begingroup%
  \makeatletter%
  \providecommand\color[2][]{%
    \errmessage{(Inkscape) Color is used for the text in Inkscape, but the package 'color.sty' is not loaded}%
    \renewcommand\color[2][]{}%
  }%
  \providecommand\transparent[1]{%
    \errmessage{(Inkscape) Transparency is used (non-zero) for the text in Inkscape, but the package 'transparent.sty' is not loaded}%
    \renewcommand\transparent[1]{}%
  }%
  \providecommand\rotatebox[2]{#2}%
  \ifx\svgwidth\undefined%
    \setlength{\unitlength}{534.5bp}%
    \ifx\svgscale\undefined%
      \relax%
    \else%
      \setlength{\unitlength}{\unitlength * \real{\svgscale}}%
    \fi%
  \else%
    \setlength{\unitlength}{\svgwidth}%
  \fi%
  \global\let\svgwidth\undefined%
  \global\let\svgscale\undefined%
  \makeatother%
  \begin{picture}(1,0.78166511)%
    \put(0,0){\includegraphics[width=\unitlength]{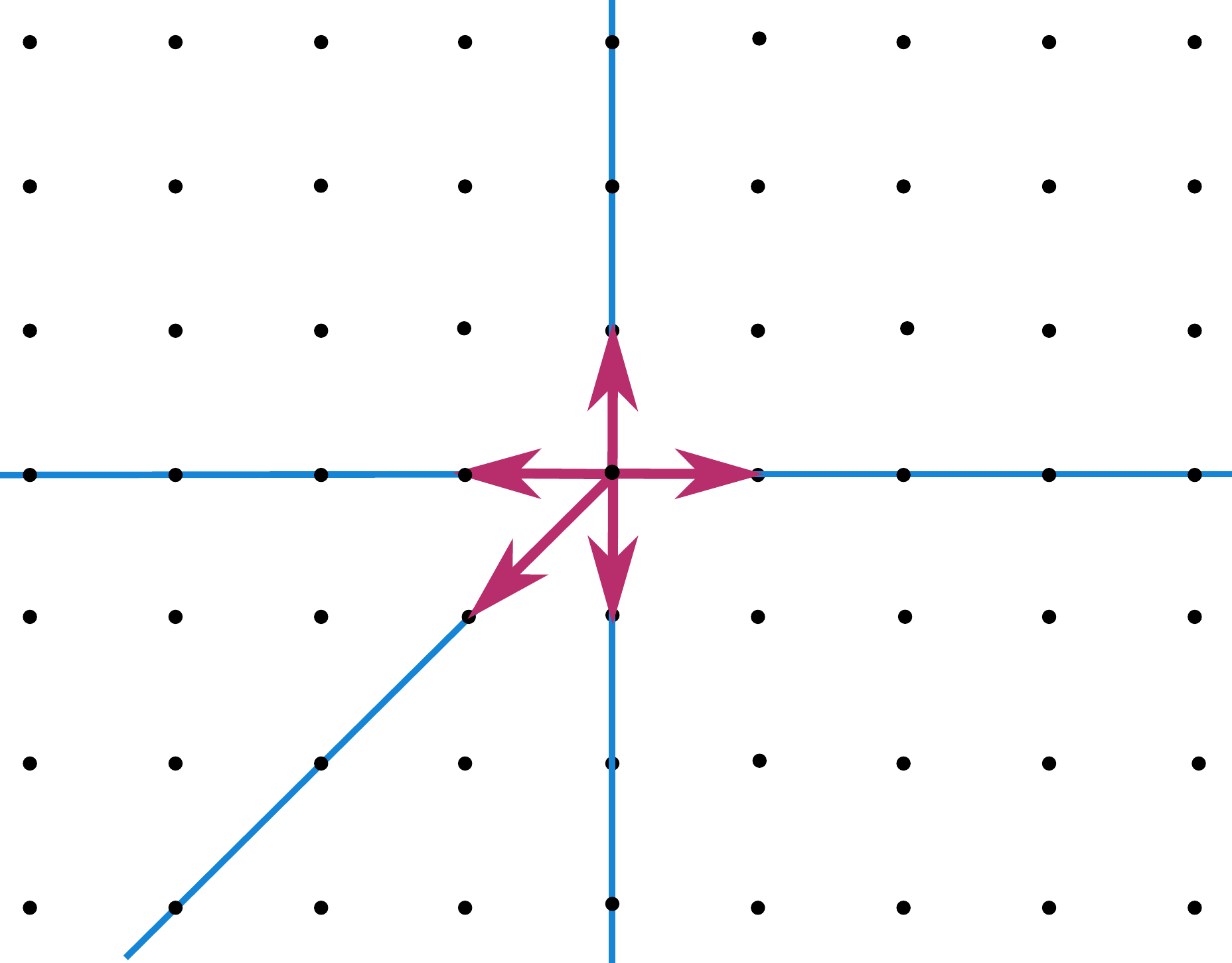}}%
    \put(0.32842943,0.30049321){\color[rgb]{0,0,0}\makebox(0,0)[lb]{\smash{{\color{mypurple} $v_4$}
}}}%
    \put(0.52089688,0.48292221){\color[rgb]{0,0,0}\makebox(0,0)[lb]{\smash{{\color{mypurple} $v_2$}
}}}%
    \put(0.58278491,0.43046843){\color[rgb]{0,0,0}\makebox(0,0)[lb]{\smash{{\color{mypurple} $v_1$}
}}}%
    \put(0.36791823,0.42953809){\color[rgb]{0,0,0}\makebox(0,0)[lb]{\smash{{\color{mypurple} $v_3$}
}}}%
    \put(0.52072245,0.30166701){\color[rgb]{0,0,0}\makebox(0,0)[lb]{\smash{{\color{mypurple} $v_5$}
}}}%
  \end{picture}%
\endgroup%

 \end{center}
We may now wonder what is preserved if each of the five rays of the above fan is rotated, so as to obtain a new fan, whose rays divide the plane into five congruent cones. The vectors of $V$ are rotated as well, into new vectors, while we may keep the same triangulation. Thus we obtain a new triangulated vector configuration $\{V',\bT\}$, whose vectors are generators of the new fan rays. However,
the $\Z$-Span of $V'$ is not a lattice but a $\Z$-module of higher rank, dense in $\R^2$. We could re-scale the vectors
of $V'$, for example into five unitary vectors. This would produce a new vector configuration $\{V'',\bT\}$. But there is no rescaling such that the $\Z$-Span of $V''$ is a lattice. In fact the new fan is nonrational, that is, there is no lattice that has non-empty intersection with each of its rays. 
\begin{center}
\def\svgwidth{0.5\textwidth}
\begingroup%
  \makeatletter%
  \providecommand\color[2][]{%
    \errmessage{(Inkscape) Color is used for the text in Inkscape, but the package 'color.sty' is not loaded}%
    \renewcommand\color[2][]{}%
  }%
  \providecommand\transparent[1]{%
    \errmessage{(Inkscape) Transparency is used (non-zero) for the text in Inkscape, but the package 'transparent.sty' is not loaded}%
    \renewcommand\transparent[1]{}%
  }%
  \providecommand\rotatebox[2]{#2}%
  \ifx\svgwidth\undefined%
    \setlength{\unitlength}{155.36447467bp}%
    \ifx\svgscale\undefined%
      \relax%
    \else%
      \setlength{\unitlength}{\unitlength * \real{\svgscale}}%
    \fi%
  \else%
    \setlength{\unitlength}{\svgwidth}%
  \fi%
  \global\let\svgwidth\undefined%
  \global\let\svgscale\undefined%
  \makeatother%
  \begin{picture}(1,0.79425097)%
    \put(0,0){\includegraphics[width=\unitlength]{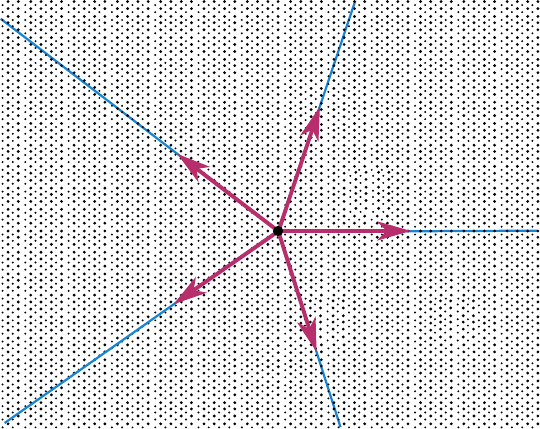}}%
  \end{picture}%
\endgroup%

 \end{center}
Hence, the notion of triangulated vector configuration still makes sense in the nonrational case, but there is no lattice.
The idea, due to Prato, is to replace the lattice with a {\em quasilattice}, that is the $\Z$-Span of a generating set 
of $\R^n$. 
Then the convex datum becomes a triple---which we will call Prato's datum---given by: a fan, a choice of rays generators, 
a choice of a quasilattice containing these generators \cite{p}. 
Notice that, with such a choice, each ray intersects the quasilattice, however, the notion of primitive vector does not make sense any longer.
As shown in the above examples, a Prato's datum can be naturally encoded in a single convex object: a triangulated vector configuration. Notice that there are many triangulated vector configurations that encode a given Prato's datum.
The theory of vector configurations, as developed in \cite{DL-R-S}, provides a precise link between triangulated vector configurations, sets of holomorphic vector fields as used in LVMB theory, fans and convex polytopes. In particular Gale duality plays an important role in connecting these different aspects. This view-point was developed in \cite{bz1}.  
The vector configurations that we consider must satisfy two further technical conditions---already considered in
\cite{MV}---they must be balanced and odd: balanced means that the sum of the vectors of $V$ is zero, moreover we require that $n-d$ is odd. 
This does not mean a loss in generality. For example, let 
$$\Big\{(\ (1,0),(0,1),(-1,1),(-1,0),(0,-1)\ ),\bT\Big\}$$
be the triangulated vector configuration considered above, which is not balanced. Consider
the new triangulated vector configuration 
$$\Big\{ (\ (1,0),(0,1),(-1,1),(-1,0),(0,-1),(1,-1),(0,0)\ ),\bT\Big\}.$$
The $\Z$-Span of the vectors is still $\Z^2$ but
$n=7$, $h=5$. The first $5$ vectors are the same as before, namely the primitive generators of the fan rays. 
The new vectors, $v_6=(1,-1)$ and $v_7=(0,0)$, will
be referred to as ghost vectors. The new triangulated vector configuration is balanced and odd.
By Gale duality, we obtain a (non unique) $\C^{2}$-action on $\C^{7}$. This action, together with $\bT$, in turn gives rise to an LVM manifold $(N,\calF)$ of complex dimension $4$. 
The leaf space $N/\calF$ is biholomorphic to $\C P^1\times\C P^1$ blown up at one point, 
associated with the rational fan drawn at page 3.
Consider now the triangulated vector configuration $\{V'',\bT\}$, whose five vectors are the five
roots of $z^5=1$ in $\C$. It is balanced and odd.
As above this yields a (non unique) LVM manifold $(N,\calF)$ of complex dimension $3$.
The leaf space $N/\calF$ is biholomorphic to the toric quasifold of complex dimension $2$ described in \cite[Example~2.9]{cx}.
In general, in \cite{bz1}, for each given odd, balanced, triangulated vector configuration $\{V,\bT\}$, not necessarily rational or polytopal, we are able to associate, via Gale duality, a (non unique) $\C^{m}$-action on $\Cn$.
This, together with $\bT$, in turn determines an LVMB manifold $(N,\calF)$. However, the complex leaf space $N/\calF$, endowed with the quotient topology, only depends on the Prato's datum encoded in the triangulated vector configuration, whilst its cohomology only depends on the combinatorial type of the fan \cite{Bat,bz1}

More precisely, let $\Delta$ be the fan associated with $\{V,\bT\}$, that is whose rays are generated by the first $h$ vectors of $V$. When the $\Z$-Span of $V$ is a lattice and the first $h$ vectors, that is, the rays generators, are primitive, the leaf space is the simplicial toric variety associated with the rational fan $\Delta$. This may be either a manifold or may have finite quotient singularities. If some of the first $h$ vectors are not primitive, the leaf space is the toric variety
associated with $\Delta$, equipped with an equivariant orbifold structure \cite{MV}. The construction of these toric orbifolds was introduced by Lerman-Tolman in \cite{LT} in the symplectic set-up; they define the notion of labeled polytope to keep track of the vectors that are integer multiples of primitive rays generators.
When the $\Z$-Span of $V$ is a quasilattice, then the leaf space is the complex toric quasifold corresponding to the
Prato's datum encoded in $\{V,\bT\}$: we are referring here to the construction of toric quasifolds, which was introduced by Prato in the symplectic set up, together with the notion of quasifold, a generalization of orbifold \cite{p}. Complex toric quasifolds were then defined in \cite{cx}. 
 
Conversely, a toric manifold, orbifold, or quasifold is always the leaf space of an LVMB manifold. This is proved in \cite{MV} in the polytopal case, for the toric manifold and orbifold cases. In \cite{bz1} we extend this result in both directions, namely to nonpolytopal and nonrational fans. This provides also an extension of complex toric quasifolds (and therefore of complex toric orbifolds) to the nonpolytopal case.

We then focus on the polytopal case. Recall that a toric variety is projective if and only if its associated fan is polytopal. In the same spirit, for LVMB manifolds, we have the following result, due to several authors: the foliation in an LVMB manifold is transversely K\"ahler if and only if the manifold is LVM. Part of the inverse implication was proved by L\oe b and Nicolau in \cite{LN} and then extended by Meersseman to the general LVM case. The direct implication was conjectured by Cupit-Foutou and Zaffran in \cite{CZ} and recently proved by Ishida in \cite{I2}. 
We recall these results and list a series of other characterizations of polytopality.
Finally, we present a variant of the Delzant construction that naturally derives from our view point and we extend to
the nonrational setting some results by \cite{MV}. 

\section{Preliminaries}

\subsection{Construction of LVMB-manifolds} \label{construction} 

We briefly describe the constructions of LVMB manifolds. Consider a configuration of points $\Lambda=(\Lambda_{1},\dots, \Lambda_{n})$ in affine space $\C^m$, that is, a finite ordered list. Repetitions are allowed but  
we assume that $\Lambda$ is not contained in a proper affine subspace. 
We define a row vector for each $j=1,\ldots,n$:
$\lr{j}:=
[\horiz\text{Re}(\Lambda_j)\horiz\ \ \horiz \text{Im}(\Lambda_j) \horiz]
\in \R^{2m}$.
A {\em basis} is a subset $\tau^*$ of $\{1,\cdots,n\}$, of cardinality $2m+1$, such that 
the interior $\mathring{C}_\alpha$ of the convex hull of
$(\lr{j})_{j\in\tau^*}$ is non empty. 
Now let us pair the configuration $\Lambda$ with a combinatorial datum, as we do when, for a given vector configuration, 
we take a triangulation. A {\em virtual chamber} $\bT^*$ of the 
configuration $\Lambda$ is a
collection of bases ${\big\{\tau^*_\alpha\big\}}_\alpha$ that satisfy Bosio's conditions \cite{Bos}, that is:
\begin{enumerate}
\item[(i)] $\mathring{C}_\alpha\cap 
  \mathring{C}_\beta  \neq\varnothing$ for every $\alpha, \beta$; 
\item[(ii)] for every $\tau^*_\alpha\in \bT^*$ and every 
$i \notin \tau^*_\alpha$,\\ 
there exists $j\in \tau^*$ such that 
$\big(\tau^*_\alpha  \setminus \{j\}\big)\cup \{i\}\in \bT^*$.
\end{enumerate}
We show in the picture an example with $n=6$ and $$\bT^*=\Big\{ 
\{135\}, \{246\}, \{136\}, \{235\}, \{145\}, \{146\}, \{236\}, \{245\} \Big\},$$ 
to facilitate visualization, we add colors and show on the right hand side slight translations of the $\mathring{C}_\alpha$'s:
\begin{center}
\def\svgwidth{0.8\textwidth}
\begingroup%
  \makeatletter%
  \providecommand\color[2][]{%
    \errmessage{(Inkscape) Color is used for the text in Inkscape, but the package 'color.sty' is not loaded}%
    \renewcommand\color[2][]{}%
  }%
  \providecommand\transparent[1]{%
    \errmessage{(Inkscape) Transparency is used (non-zero) for the text in Inkscape, but the package 'transparent.sty' is not loaded}%
    \renewcommand\transparent[1]{}%
  }%
  \providecommand\rotatebox[2]{#2}%
  \ifx\svgwidth\undefined%
    \setlength{\unitlength}{368.8217041bp}%
    \ifx\svgscale\undefined%
      \relax%
    \else%
      \setlength{\unitlength}{\unitlength * \real{\svgscale}}%
    \fi%
  \else%
    \setlength{\unitlength}{\svgwidth}%
  \fi%
  \global\let\svgwidth\undefined%
  \global\let\svgscale\undefined%
  \makeatother%
  \begin{picture}(1,0.33084238)%
    \put(0,0){\includegraphics[width=\unitlength]{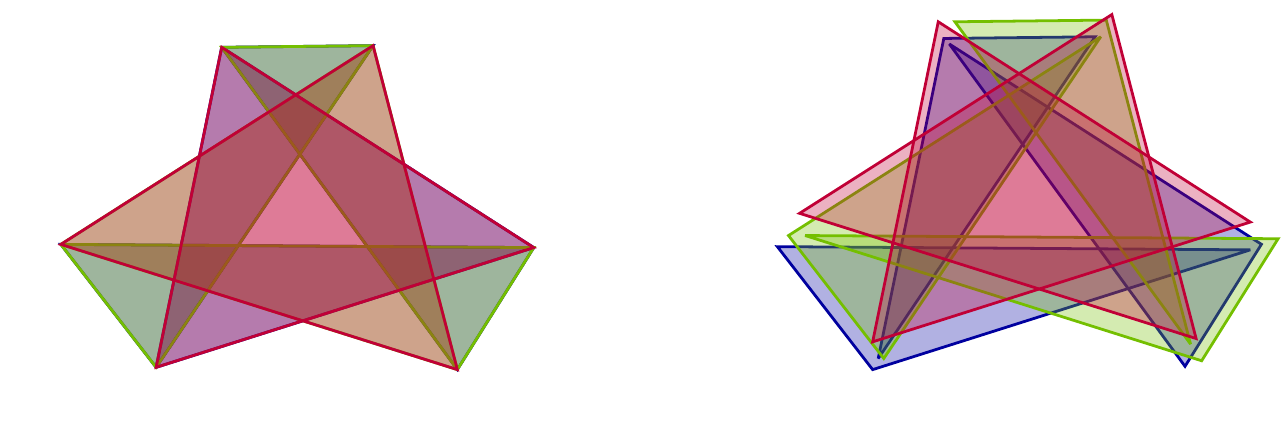}}%
    \put(0.11069696,0.00613864){\color[rgb]{0,0,0}\makebox(0,0)[lb]{\smash{$\Lambda_1$}}}%
    \put(0.35363277,0.00613864){\color[rgb]{0,0,0}\makebox(0,0)[lb]{\smash{$\Lambda_2$}}}%
    \put(0.43165335,0.14029773){\color[rgb]{0,0,0}\makebox(0,0)[lb]{\smash{$\Lambda_3$}}}%
    \put(-0.0021606,0.14029773){\color[rgb]{0,0,0}\makebox(0,0)[lb]{\smash{$\Lambda_6$}}}%
    \put(0.28815014,0.31106656){\color[rgb]{0,0,0}\makebox(0,0)[lb]{\smash{$\Lambda_4$}}}%
    \put(0.15800595,0.31106656){\color[rgb]{0,0,0}\makebox(0,0)[lb]{\smash{$\Lambda_5$}}}%
  \end{picture}%
\endgroup%

 \end{center} 
In general, $\bigcap_{\alpha} \mathring{C}_\alpha =\varnothing$.  We say that a virtual chamber is a {\em chamber} when all 
$\mathring{C}_\alpha$'s do have a common intersection. 
Example with $n=6$ and $\bT^*=\Big\{ 
\{124\}, \{134\},\{135\},\{136\},\{235\},\{236\},\{245\},\{246\} \Big\}$:
\begin{center}
\def\svgwidth{0.8\textwidth}
\begingroup%
  \makeatletter%
  \providecommand\color[2][]{%
    \errmessage{(Inkscape) Color is used for the text in Inkscape, but the package 'color.sty' is not loaded}%
    \renewcommand\color[2][]{}%
  }%
  \providecommand\transparent[1]{%
    \errmessage{(Inkscape) Transparency is used (non-zero) for the text in Inkscape, but the package 'transparent.sty' is not loaded}%
    \renewcommand\transparent[1]{}%
  }%
  \providecommand\rotatebox[2]{#2}%
  \ifx\svgwidth\undefined%
    \setlength{\unitlength}{369.1664641bp}%
    \ifx\svgscale\undefined%
      \relax%
    \else%
      \setlength{\unitlength}{\unitlength * \real{\svgscale}}%
    \fi%
  \else%
    \setlength{\unitlength}{\svgwidth}%
  \fi%
  \global\let\svgwidth\undefined%
  \global\let\svgscale\undefined%
  \makeatother%
  \begin{picture}(1,0.32619939)%
    \put(0,0){\includegraphics[width=\unitlength]{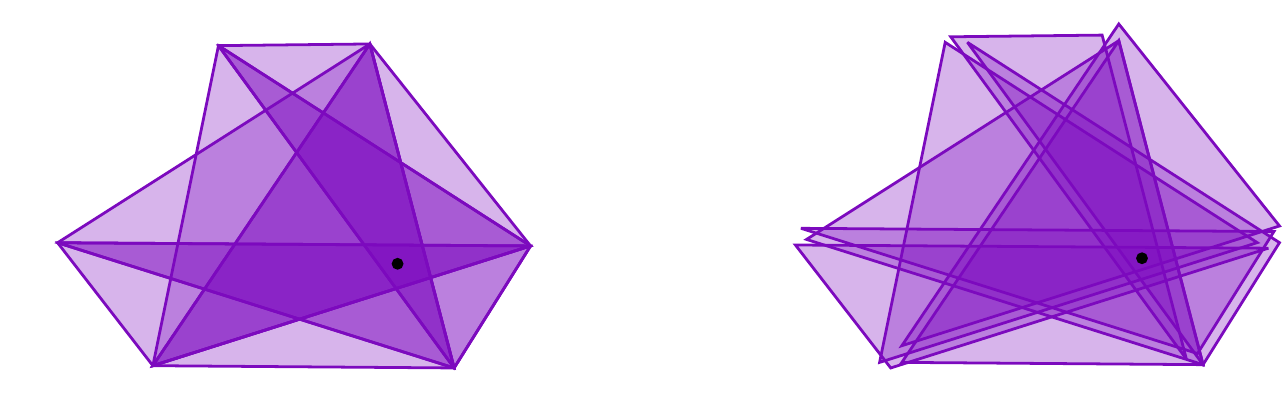}}%
    \put(0.11059358,0.0061329){\color[rgb]{0,0,0}\makebox(0,0)[lb]{\smash{$\Lambda_1$}}}%
    \put(0.35330252,0.0061329){\color[rgb]{0,0,0}\makebox(0,0)[lb]{\smash{$\Lambda_2$}}}%
    \put(0.42258206,0.13583256){\color[rgb]{0,0,0}\makebox(0,0)[lb]{\smash{$\Lambda_3$}}}%
    \put(-0.00215858,0.13583256){\color[rgb]{0,0,0}\makebox(0,0)[lb]{\smash{$\Lambda_6$}}}%
    \put(0.28788104,0.30644204){\color[rgb]{0,0,0}\makebox(0,0)[lb]{\smash{$\Lambda_4$}}}%
    \put(0.15785839,0.30644204){\color[rgb]{0,0,0}\makebox(0,0)[lb]{\smash{$\Lambda_5$}}}%
  \end{picture}%
\endgroup%

 \end{center}
We'll see below that the special case when a virtual chamber is a chamber corresponds in the toric context to a fan being polytopal (or a toric variety being projective), and in the context of LVMB theory to a manifold being transversely K\"ahler.

An {\em LVMB datum} $\{ \Lambda, \bT^{*} \}$ is a configuration 
$\Lambda=(\Lambda_{1},\dots, \Lambda_{n})$ in $\C^{m}$, with $n\geq 2m+1$, 
together with a choice of a virtual chamber $\bT^*$. 
We now show how to define a compact complex manifold $N$ from an LVMB datum. 
From the virtual chamber $\bT^*$ we define a subset $U(\bT^*)$ of $\C\P^{n-1}$ as follows: 
for each $\tau^*\in \bT^*$, define $U_{\tau^*}
:=\setofst{[z_1:\dots:z_n]\in \C\P^{n-1}}{ \forall j\in \tau^*, z_j\neq 0}$ and take 
$$U(\bT^*):=\bigcup_{\tau^*\in \bT^{*}} U_{\tau^*}.$$
We write $\Lambda$ as the matrix
 $$
\Lambda= \begin{bmatrix}
\horiz \Lambda_1 \horiz \\
   \vdots                 \\
\horiz \Lambda_{n} \horiz \\
\end{bmatrix}\in\C^{n\times m }
$$
and define the subspace $\h\subset\C^n$ as the span of the $m$ 
{\em columns} of $\Lambda$.

By Bosio's condition (i), $\h$ has dimension $m$. 
Let {$\exp\colon\C^n\rightarrow(\C^*)^n$. The {action by
$\exp\h\subset(\C^*)^n$ on $\C^n$, induced by the natural $(\C^*)^n$-action, is a $\C^m$--action that commutes} with the diagonal $\C^*$--action. Then
\cite{Bos} the group  $\exp\h$}
acts freely and properly on $U(\bT^*)$, so we can define
$$N:=U(\bT^*)/\exp\h.$$
Manifolds arising from this construction are called {\em LVMB manifolds} \cite{ldm,M, Bos}. {They have a very rich
geometry.}

Every LVMB manifold $N$ admits a noteworthy smooth holomorphic foliation $\calF$, which has been investigated by several authors \cite{LN,M, MV, CZ, B, Tam, PU, U, PUV}. Following \cite{U}, we can simply describe the leaves of $\calF$ as the orbits of the action of $\exp \overline{\h}$ 
on $N$, which descends from the action on $U(\bT^*)$. {There is also on $N$ an induced action of the abelian complex group given by the quotient of $(\C^*)^n/\exp(\h)$ by $\text{diag}(\C^*)$. This group has the same complex dimension as $N$, $n-1-m$, and has a dense open orbit in $N$. It also contains a real group isomorphic to $(S^1)^n/\text{diag}(S^1)$. 
With respect to any metric that makes this real group act by isometries, the foliation $\calF$ is Riemannian \cite[Sect.~2.3.2]{bz1}.

We thus have the following
\begin{thm}\label{the manifold N}
An LVMB datum $\{ \Lambda, \bT^{*} \}$ determines a compact complex manifold $N$ of dimension $n-1-m$. This manifold is endowed
with a holomorphic Riemannian foliation $\calF$. Moreover there is an abelian complex group of same dimension acting on $N$
with a dense open orbit. In the limiting case $n=2m+1$, the manifold $N$ is a compact complex torus, but whenever $n>2m+1$, $N$ is not K\"ahler. 
\end{thm}
}

Another entry point of LVMB theory is that of moment-angle manifolds. This point of view was developed by Tambour and Panov-Ustinovsky \cite{Tam, PU}. The latter authors and Battisti \cite{B} are precursors to the toric methods established by Ishida \cite{I2}, who also proved a remarkable group-theoretic characterization of LVMB manifolds as a large subclass of the class of complex manifolds that admit a so-called {\em maximal torus action}. Namely, the action of a compact torus $G$ on a manifold $M$ is called {\em maximal} when there exists $x\in M$ such that 
$\dim_{\R} G + \dim_{\R} G_{x}=\dim_{\R} M$, where $G_{x}$ is the isotropy subgroup at $x$. For LVMB manifolds, the torus $G$ is the above-mentioned real group isomorphic to $(S^1)^n/\text{diag}(S^1)$.
 In \cite{U} Ustinovsky shows that the family of complex manifolds presented in \cite{PU} coincides with Ishida's complex manifolds admitting a maximal torus action.

\subsection{Duality LVMB-data $\leftrightarrow$ fans and (quasi)lattices} \label{link to toric}
In this section we illustrate how, via Gale duality, we are able to pass from a set of data to the other.
Let $W=(w_{1},\dots,w_{n})$ denote a vector configuration in some real, finite-dimensional vector space. 
We define the set of relations of $W$ by 
$$\text{Rel}(W):=\setofst{(\alpha_{1},\dots,\alpha_{n})}{\sum \alpha_{j} w_{j}=0}\subset \R^{n}.$$  
We say that a vector configuration $W$ is {\em balanced} when $\sum w_{j}=0$; we say that it is {\em graded} when there exists an affine hyperplane not passing through the origin that contains all vectors of $W$. A graded configuration of vectors is a useful way to represent and manipulate a configuration of {\em points}, i.e. elements of an affine space.

Now consider an LVMB datum $\{ \Lambda, \bT^{*} \}$. Let $A$ be a $2m$-dimensional affine subspace of $\R^{2m+1}$ not passing through the origin. 
Send the configuration $\Lambda^\R=(\lr{1},\dots,\lr{n})$, defined in Section~\ref{construction}, in $A$ via an 
affine real isomorphism from $\R^{2m}$ to $A$. 
This {determines}
a graded vector configuration $W(\Lambda)$ in $\R^{2m+1}$. Consider a $(n-2m-1)\times n$ matrix whose rows are a basis of 
$\text{Rel}(W(\Lambda))$.
The $n$ columns of this matrix form a vector configuration in $\R^{n-2m-1}$ that we denote $V$. This configuration is said to be in Gale duality with the configuration $\Lambda^\R$. The configuration $V$ is balanced, and also odd since $n-(n-2m-1)$ is odd.

The collection $\bT$ of all the complements (in $\{1,\dots,n\}$) of elements of $\bT^{*}$ is 
a {\em triangulation} of $V$ (cf.  \cite{DL-R-S}). 
Roughly speaking, a triangulation of $V$ is the simplicial complex determined by a simplicial fan whose ray generators form a subset of $V$.

Therefore, we have obtained, from $\{ \Lambda, \bT^{*} \}$, a
triangulated vector configuration $\{V,\bT\}$, which is 
 unique up to an ambient real linear automorphism--{for our purposes we can consider $\{V,\bT\}$ uniquely determined
 by $\{ \Lambda, \bT^{*} \}$.} 
Elements of $\bT$ are called simplices.
{As in \cite{M} each index in the intersection of all of the bases of $\bT^*$ is
an {\em indispensable index} of $\bT^*$. We denote by $k$ the number of such indices and we will always assume that $j\in \{1,\cdots,n\}$ is indispensable if and only if $j> n-k$.
 Any indispensable index $i$ of $\bT^{*}$ will not appear in any simplex of $\bT$; it is called a {\em ghost index}  of $\bT$, and $v_{i}$ is called a {\em ghost vector}. By properties of Gale duality \cite[4.1.38 (iv)]{DL-R-S}, for any $\tau\in \bT$ the vectors $(v_{j})_{j\in\tau}$ are linearly independent, so they generate a simplicial cone denoted $\text{cone}(\tau)$.  One can check that Bosio's conditions 
are equivalent to the fact that the collection of these cones for all $\tau\in \bT$ determines a complete 
simplicial fan $\Delta$ \cite[Prop. 2.1]{bz1} ; this fan is not necessarily rational. The non necessarily closed additive subgroup $Q$ of $\R^{n-2m-1}$ generated by 
$v_{1},\dots,v_{n}$ is called a {\em quasilattice }(\cf \cite{p}). We always see $Q$ as embedded, i.e. $Q$ is a shorthand for the pair $(\R^{n-2m-1},Q)$. 

In the introduction we have described how a toric datum, and, more precisely, a Prato's datum, can be encoded in a triangulated vector configuration. In more detail, consider the triple $(\Delta,\{v_1,\ldots,v_h\},Q)$ where $\Delta\subset\R^d$ is a complete fan, $\{v_1,\ldots,v_h\}$ is a set of fan rays generators and $Q$ is a quasilattice in $\R^d$ containing these generators. Then a corresponding triangulated vector configuration, balanced and odd, is a pair 
$\{V,\bT\}$ with $V=(v_1,\ldots,v_n)$ a vector configuration in $\R^d$ such that the first $h$ vectors coincide
with the given generators above, $Q=\text{Span}_{\Z}\{v_1,\ldots,v_n\}$, $\sum_i v_i=0$ and $n-d=2m+1$, with $m\in\N$.
On the other hand $\bT$, as we have seen, corresponds to the fan cones and it is a combinatorial datum. 
Remark that we may have to add vectors to the set $\{v_1,\ldots,v_h\}$ in order to have all the above conditions satisfied.
Notice that, in the classical toric setting, $Q$ is a fixed lattice and the generators are taken to be primitive in $Q$. 
However, here we can take non primitive generators , and, in fact, any set of generators is allowed, 
as long as the quasilattice $Q$ is chosen consistently. 

Once a balanced and odd triangulated vector configuration on $\{V,\bT\}$ is given, Gale duality yields a (non unique) LVMB datum as follows: 

Consider a $(n-d=2m+1)\times n$ matrix whose rows are a basis of 
$\text{Rel}(V)$. 
The $n$ columns of this matrix form a vector configuration in $\R^{2m+1}$ that we denote $\Lambda^\R$, unique up to an ambient real linear automorphism. Reversing the construction of Section~\ref{construction} yields the configuration 
$\Lambda$ in $\C^{m}$, which is not unique up to an ambient complex affine automorphism {(see \cite[Section~2.2.2]{bz1} for details)}.  
The complements of bases in the triangulation $\bT$ form a virtual chamber $\bT^{*}$, so we obtain an LVMB datum $\{\Lambda, \bT^{*}\}$. Due to the non uniquess of $\Lambda$, we actually obtain from 
$\{V,\bT\}$ a family of (closely related) LVMB manifolds. We propose to call two members of such a family {\em virtually biholomorphic}. 
{The questions posed by the non-uniqueness of the Gale duality correspondence are treated in our paper \cite{bz2}, in preparation.

To complete this preliminary section let us recall how to measure the nonrationality of a vector configuration $V$. Consider the space of linear relations  $\Rel(V)\subset \R^{n}$, of dimension $n-d$. We say that a real subspace of $\R^{n}$ is {\em rational} when it admits 
a real basis of vectors in $\Q^{n}$ (equivalently, $\Z^{n}$). 
We define $a(V)$ as the dimension of the largest rational space contained 
in $\Rel(V)$, and $b(V)$ as the dimension of the smallest rational space 
containing $\Rel(V)$ (our number $a(V)$ is closely related to the number $a$ defined in \cite{M}). Then $0\leq a(V) \leq n-d \leq b(V) \leq n$.

The configuration is called {\em rational} when $\Rel(V)$ is rational or, equivalently,         
$a(V)= n-d$, or $b(V)= n-d$. Otherwise $2+ a(V)\leq b(V)$, and all such values are possible. 

Notice that $\text{Span}_{\Z}(V)$ is an honest lattice if and only if the configuration $V$ is rational.
These numbers determine the topology of the generic leaves of the foliation $\calF$ and their closures \cite[Sect. 2.3]{bz1}, for example generic leaves are homeomorphic to 
${(S^{1})}^{a(V)-1}\times \R^{2m-a(V)+1}$.
If the configuration is rational, that is $a(V)=b(V)=2m+1$, all leaves are closed 
(see \cite{MV} for a full study of this case).
On the other hand there are nonrational configurations $V$ such that $a(V)=1$; in these cases the generic 
leaf is $\C^m$.

\section{Toric and LVMB geometry} \label{Toric and LVMB}
In this section we present, in form of synthetic statements, some of the results that we proved in \cite{bz1}, adding some new features. For each result we outline the argument of the proof.
\begin{thm}\label{leaf space} Let $\{\Lambda,\bT^*\}$ be an LVMB-datum and
let $(N,\calF)$ be the corresponding LVMB manifold. Let $\{V,\bT\}$ be the triangulated vector configuration 
Gale dual to $\{\Lambda,\bT^*\}$. Let $(\Delta,Q,\{v_1,\ldots,v_h\})$ be the fan, quasilattice and ray generators encoded in $\{V,\bT\}$ and let $X$ be the toric quasifold associated to this triple.
Then the leaf space $N/\calF$ is biholomorphic to $X$.
\end{thm}
First of all let us recall how we obtain $X$. From the Audin-Cox-Delzant construction \cite{A,cox,D} and its nonrational complex generalization \cite{cx}, it is known that to $(\Delta,Q,\{v_1,\ldots,v_h\})$ there corresponds a geometric quotient $X=U'(\Delta)/G$, where $U'(\Delta)$ is an open subset of $\C^h$ that 
depends on the combinatorics of $\Delta$, and $\text{N}_\C$ is a complex subgroup of $(\C^*)^h$ that depends on $Q$ and on the 
vectors $v_1,\ldots,v_h$. If the configuration is rational (resp. nonrational), then $X$ is 
a complex manifold or a complex orbifold (resp. a non Hausdorff complex quasifold) of dimension $d$,
acted on holomorphically by the torus (resp. {\em quasitorus}) $\C^d/Q$ (\cf \cite{A,cox,p,cx}); 
the construction in \cite[Thm~2.2]{cx} can be adapted to the nonpolytopal case.
Quasifolds generalize orbifolds: the local model is a quotient of a manifold by the smooth action of a finite or countable group, non free on a closed subset of topological codimension at least $2$ (\cite{p}, see also \cite{kite}). 
Let us describe more precisely how we construct $X$, and its relation with $N$.
Let 
\begin{equation}\label{uhat}
\hat{U}(\bT^*)=:\bigcup_{\tau^*\in \bT^{*}}\setofst{(z_1,\dots,z_n)\in \C^{n}}{ \forall j\in \tau^*, z_j\neq 0}.
\end{equation}
Then 
$U(\bT^*)$ is just the projectivization of $\hat{U}(\bT^*)$. Consider
$\C^n=\C^h\times\C^k$, then $U'(\Delta)$ is exactly the projection
onto the first factor of $\hat{U}(\bT^*)$. The manifold $N$ is given by the quotient $U(\bT^*)/\exp(\h)$. In order to pass from $\{\Lambda,\bT^*\}$ to $\{V,\bT\}$, let us identify $\R^{2m}$ with the affine space $A=\{(1,x_1,\ldots,x_{2m})\}\subset\R^{2m+1}$ (cf. Section~\ref{link to toric}). Then} the leaf space $N/\exp(\overline{\h})$ can be naturally identified with 
$$\hat{U}(\bT^*)/\exp(Rel(V)_{\C}).$$
On the other hand $X$ is the orbit space $U'(\Delta)/\text{N}_{\C}$, endowed with the quotient topology. Here  
\begin{equation}\label{enne}\text{N}=\exp\Big(\{\underline{a}\in\R^h\;|\;\sum_{i=1}^h a_iv_i\in Q\}\Big)\end{equation} and
\begin{equation}\label{enneci}\text{N}_{\C}=\exp\Big(\{\underline{a}\in\C^h\;|\;\sum_{i=1}^h a_iv_i\in Q\}\Big)\end{equation}
{are subgroups of $(S^1)^h$ and $(\C^*)^h$ respectively} (see \cite{p,cx} for details).
When $h=n$, we have $\text{Span}_{\Z}\{v_1,\ldots,v_h\}=Q$, therefore $\exp(Rel(V)_\C)=\text{N}_{\C}$. Moreover, $\hat{U}(\bT^*)=U'(\Delta)$. Thus the leaf space $N/\calF$
and $X$ are naturally identified. When $h<n$, we can define the map $f\colon N/\calF\rightarrow X$ as follows:
let $[z_1,\dots,z_n]\in\hat{U}(\bT^*)/(\exp(Rel(V)_\C)$. {Notice that, by (\ref{uhat}), $z_j\neq0$ for all $j=h+1,\ldots,n$.} Because of the properties of $(V,\bT)$, there exists 
$\underline{b}\in Rel(V)_{\C}$ such that $e^{2\pi i b_j}\cdot z_j=1$, for all $j=h+1,\ldots,n$.
The mapping 
$$f([z_1,\dots,z_h,z_{h+1},\dots,z_n])=[e^{2\pi i b_1}z_1,\dots,e^{2\pi i b_h}z_h]$$
is well defined and
identifies the leaf space $N/\calF$ and $X$ as complex quotients.
Let us check these properties of $f$. Let $\underline{b},\underline{b}'\in Rel(V)_{\C}$ such that 
$e^{2\pi i b_j}\cdot z_j=e^{2\pi i b'_j}\cdot z_j=1$, for all $j=h+1,\ldots,n$. Then 
$\sum_{j=1}^h(b_j-b'_j)v_j\in\text{Span}_\Z\{v_{h+1},\dots,v_n\}\subset Q$. 
Therefore $(e^{2\pi i (b_1-b'_1)},\dots,e^{2\pi i (b_h-b'_h)})\in \text{N}_\C$, thus the map $f$ does not depend on the choice of $\underline{b}$. Similarly it can be checked that $f$ does not depend on the class representative $\vz$ and that it is injective. Surjectivity is immediate. We can construct complex atlases for the two quotients $\hat{U}(\bT^*)/\exp(Rel(V)_\C)$ and
$U'(\Delta)/\text{N}_\C$ very similarly (cf. \cite[Theorem~2.1]{cx} and the proof of 
\cite[Lemma~3.2]{bz1}), namely using holomorphic slices for the action of the respective groups given by
$(z_1,\ldots,z_d)\rightarrow (z_1,\ldots,z_d,\underbrace{1,\ldots,1}_{n-d})$ and 
$(z_1,\ldots,z_d)\rightarrow (z_1,\ldots,z_d,\underbrace{1,\ldots,1}_{h-d})$.
The map $f$, locally, is just the identity. We deduce that, in particular,
the complex structure induced by $(N,\calF)$ on the leaf space depends only on $(\Delta,Q,\{v_1,\ldots,v_h\})$.
{In particular it does not depend on the choice of ghost vectors, nor on the choice of the Gale dual LVMB datum. 
This is consistent with \cite[Theorem~G,(iii)]{MV}. 
The (quasi)torus $$\C^d/Q\simeq\C^h/\text{N}_\C\simeq\C^n/\exp(Rel(V)_\C)$$ acts on $X$ holomorphically with a dense open orbit.
\begin{rem} If $Q$ is a lattice then $X$ is a toric orbifold, if $Q$ is a lattice and the $v_i$ are primitive then
$X$ is the toric variety associated to $\Delta$, with no additional orbifold structure.
\end{rem}
\begin{rem}
Remark that the action of the group $\text{N}$ does induce a holomorphic foliation on $U'(\Delta)$.  
However, since $\text{N}$ is, in general, for rational and nonrational configurations, not connected,
the leaf space is {\em not} $X$ (see for a nonrational example \cite{rhombus} and for a rational example \cite[Example~2.4.2]{bz1}). Following \cite{MV}, this problem is overcome in our construction by "increasing the dimension" (cf. \cite[Example~2.4.2]{bz1}).

We will see this same phenomenon in the symplectic setting with the variant of the Delzant construction introduced in Section~\ref{delzant}.
\end{rem}
\begin{thm} Let $(\Delta,Q, \{v_1,\ldots,v_h\})$ be a complete fan, a quasilattice and ray generators in $Q$. 
Let $X$ be the toric quasifold associated to this triple. Then there exists an  LVMB manifold $(N,\calF)$ whose leaf space is biholomorphic to $X$; there are infinitely many such LVMB manifolds. 
\end{thm}
By Section~\ref{link to toric} it is sufficient to consider a triangulated vector configuration $\{V,\bT\}$ encoding 
$(\Delta, Q,\{v_1,\ldots,v_h\})$. We then take any Gale dual LVMB--datum $\{\Lambda,\bT^*\}$. This determines an LVMB manifold $N$, whose leaf space $N/\calF$ can be identified with $X$ by Theorem~\ref{leaf space}.
The next two corollaries are proved in the polytopal case in \cite[Theorem G]{MV}; the first is proved, in this generality, in \cite{CZ}.
\begin{crl} Let $(\Delta,L,\{v_1,\ldots,v_h\})$ be a complete fan, lattice and primitive ray generators in $L$ and 
let $X$ be the toric variety determined by this triple. Then there exists an infinite family of LVMB manifold $(N,\calF)$ whose leaf space is biholomorphic to $X$.
\end{crl}
\begin{crl}
 Let $(\Delta,L,\{v_1,\ldots,v_h\})$ be a complete fan, lattice and ray generators in $L$ and 
let $X$ be the toric variety obtained from this triple, with the orbifold structure determined by the ray generators multiplicities. Then there exists an infinite family of LVMB manifolds $(N,\calF)$ whose leaf space is biholomorphic to $X$, with its orbifold structure.
\end{crl}

In conclusion, LVMB manifolds, through Gale duality, model rational and nonrational complete simplicial varieties, not necessarily projective, and avoid all singularities \cite{bz1}. 
We note also the interesting announcement of Katzarkov--Lupercio--Meersseman--Verjovsky \cite{KLMV} that the nonrational LVM manifolds can also serve as a model for a notion of noncommutative (projective, simplicial) toric varieties.
\section{The polytopal case}\label{polytopal}
{In this section we specialize to polytopal fans:
we say that a fan $\Delta$ is polytopal when it is the {\em dual} (or {\em normal}) fan to some polytope $P$. Loosely speaking, this happens when there exists a polytope $P$ such that the fan rays are normal to the polytope facets, and the other cones of $\Delta$ are related to the faces of $P$ by an inclusion-reversing bijection. This implies
important properties for the corresponding spaces. For example it is well known that
a toric variety is projective if and only if its fan is polytopal. 
An analogous result, whose proof was completed by Ishida in 2015, holds at the level of LVMB manifolds. First we need
to recall a definition: a complex foliated manifold $(N,\calF)$ is {\em transversely K\"ahler} when there exists a closed $(1,1)$ form $\omega$ on $N$ such that {$\omega_x(X,JX)\geq0$} for all $x\in X$ and for all $X\in T_xX$, with equality if and only if $X$ is tangent to $\calF$. We then have:}
\begin{thm}\label{transverse}
Let $(\Lambda_1,\ldots,\Lambda_n)$ be an affine point configuration in $\C^m$, together with the choice of a virtual chamber $\bT^*$, satisfying Bosio's conditions. Let $\Delta$ be the associated fan. Then the LVMB manifold $(N,\calF)$ 
corresponding to $(\Lambda_1,\ldots,\Lambda_n)$ is transversely K\"ahler if and only if $\Delta$ is polytopal.
\end{thm}
Proof of the inverse implication by  \cite{LN} (for $m=1$) and \cite{M} (for $m\geq 2$); proof of the direct implication by \cite{CZ} under condition (K). The general case is due to \cite[Theorem~5.7]{I2}.

Polytopality has several further characterizations, depending on the view-point we want to stress. Here is a list:
\begin{enumerate}
\item[$\phantom{\Leftrightarrow}$\quad  1.]  the fan $\Delta$ is {\em polytopal};
\item[$\Leftrightarrow$\quad  2.] the foliation $\calF$ is transversely K\"ahler;
\item[$\Leftrightarrow$\quad 3.] The triangulation is {\em regular}
\item[$\Leftrightarrow$\quad  4.] There exists a {\em height function} on $V$ that induces $\bT$
\item[$\Leftrightarrow$\quad 5.] The virtual chamber defines a nonempty {\em chamber}, i.e.,  
$\bigcap_{\alpha} \mathring{C}_\alpha \neq\varnothing$ 
\item[$\Leftrightarrow$\quad 6.] 
There exists $\nu\in \R^{2m}$  such that $\forall \tau\subset \{1\dots n\}, \tau\in \bT$ if and only if $\nu$ is 
in the interior of the convex hull of  $\setofst{\Lambda^{\R}_{j}}{j\in\tau^{c}}$
\end{enumerate}
Note that 4. is the explicit definition of 3. We refer to \cite{DL-R-S} for details on triangulations and the meaning of regularity.

{Let us mention, in conclusion, some result about the topology of the leaf space. The Betti numbers of a complex toric
quasifold are computed in \cite[3.3]{Bat}, under the polytopality assumption: odd Betti numbers are zero, even Betti numbers give the $h$-vector of the corresponding fan. This result is proved for the larger class of simplicial shellable fans in \cite[Theorem~3.1]{bz1}, where the basic cohomology of the foliated manifold $(N,\calF)$ yields
the cohomology of the leaf space. We also prove, by adding back polytopality, that the  basic cohomology is generated in degree two. Furthemore El Kacimi's theorem \cite[3.4.7]{EK} applies and gives the Hard Lefschetz theorem for the basic cohomology of LVM manifolds.}
\section{A variant of the Delzant construction} \label{delzant}
The construction introduced by Delzant in \cite{D} allows to explicitly construct, via symplectic reduction, a symplectic toric manifold from a Delzant polytope. This is a convex polytope of full dimension in $(\R^d)^*$ that is simple, rational and satisfies a certain integrality condition. Roughly speaking
a Delzant polytope $P$, at each vertex, looks like an orthant at the origin. Simple means that, for each vertex $\nu$, there are exactly $d$ facets of $P$ meeting at $\nu$.
The fan $\Delta$ normal to $P$ is the fan $\Delta\subset\R^d$ whose rays are normal to the polytope facets and have inward pointing directions. The polytope $P$ is simple if and only if the fan $\Delta$ is simplicial.
The polytope $P$ is rational if there exists a lattice $L$ in $\R^d$ with respect to which $\Delta$ is rational.
Duality yields a bijective correspondence between the polytope vertices and the fan maximal cones.
The Delzant's integrality condition is satisfied if, for each vertex $\nu$, the primitive ray generators of its corresponding maximal cone give a basis of $L$. On the other hand a symplectic toric manifold 
$(M,\omega)$ is a compact, connected, symplectic manifold equipped with the effective Hamiltonian action of a torus $T$ such that $\dim M=2\dim T$. By torus in this section we mean a compact torus, that is a torus isomorphic to $(S^1)^r$, for some $r\in\N_{>0}$. 
The convexity theorem of Atiyah \cite{a} and Guillemin--Stenberg \cite{gs} asserts that, 
if $(M,\omega)$ is a compact, connected, symplectic manifold, endowed with the Hamiltonian action of a torus $T$, with group lattice $L$ and Lie algebra $\t=L\otimes_{\Z}\R$, then the image of the corresponding moment mapping $\Phi$ is a rational convex polytope in $\t^*$, called moment polytope. A very important application of the convexity theorem is the Delzant's theorem, that completely classifies symplectic toric manifolds: the moment polytope of a symplectic toric manifold $(M,T,\omega,\Phi)$ is a Delzant polytope $P$; in turn, $P$ uniquely determines $(M,T,\omega,\Phi)$
up to equivariant symplectomorphims. A key point in the proof of this theorem is the above-mentioned procedure, that allows to construct, from a given Delzant polytope, a symplectic toric manifold. This has proved to be an extremely fruitful tool, in symplectic and contact geometry, with a great variety of applications.

We present here a simple variant of this construction, that, as we will see, also applies to the generalizations of the Delzant procedure introduced by Lerman-Tolman \cite{LT} and Prato \cite{p}. Furthemore, our variant is strictly related to the LVM manifolds described in the previous section. 
But let us recall first the classical Delzant construction.
Let $L$ be a lattice of rank $d$ in $\R^d$. A convex polytope $P$ in $(\R^d)^*$ can be always written as intersection of closed half spaces. When these closed half spaces are in bijective correspondence with the polytope facets, this intersection is said to be minimal:
$$P=\cap_{j=1}^{h}\{\mu\in (\R^d)^*\;|\;\langle\mu,v_j\rangle\geq l_j\}.$$
Here the $v_i$ are taken to be primitive and inward pointing, the $l_i$'s are real coefficients determined by the $v_i$. Denote $V=(v_1,\ldots,v_h)$.
Let  $\text{N}$ be the subtorus of $(S^1)^h$ defined in (\ref{enne}). The integrality condition implies that $\text{N}=\exp(Rel(V))$. The induced action
of this group on $\C^h$ is Hamiltonian with respect to the standard K\"ahler form of $\C^h$. Let 
$\Psi\colon \C^h\rightarrow (Rel(V))^*$ be the corresponding moment mapping 
(the choice of the constant is such that $\Psi(0)=i^*(\sum_{j=1}^hl_je_j^*)$, where $i\colon Rel(V)\hookrightarrow \R^h$ is the inclusion).
Then the reduced space $\Psi^{-1}(0)/\text{N}$ is the toric symplectic manifold $(M_P,\R^d/L,\omega,\Phi)$ of dimension $d$, with moment polytope $P$. Notice that $\R^d/L\simeq(S^1)^h/\text{N}$.

The question that we want to pose is: what happens if we consider an intersection of half spaces that is not minimal? And what is the relation of this set-up with LVM manifolds and the previous sections?
Let us start with the following
\begin{prop}\label{sopra} Consider $n\geq h$. Let $L$ be a lattice in $\R^d$ and let $P$ be the Delzant polytope in $(\R^d)^*$
defined by a not necessarily minimal intersection 
$$P=\cap_{j=1}^{n}\{\mu\in (\R^d)^*\;|\;\langle\mu,v_j\rangle\geq l_j\},$$
where
\begin{itemize}
\item $v_1\ldots,v_h$ are primitive in $L$ and generate the $h$ 
rays of the normal fan to $P$;
\end{itemize}
and, if $n>h$,
\begin{itemize}
\item for each $j=h+1,\ldots,n$, $v_j\in L$ 
\item $P\subset\{\mu\in (\R^d)^*\;|\;\langle\mu,v_j\rangle>l_j\},$ with $j=h+1,\ldots,n$
\end{itemize}
Consider $V=(v_1,\ldots,v_n)$ and take the subgroup
$\text{N}=\exp(Rel(V))\subset(S^1)^n$. Let $\Psi$ be the moment mapping with respect to
the induced action of $\text{N}$ on $\C^{n}$. Then the reduced
space $M=\Psi^{-1}(0)/\text{N}$ is a compact symplectic manifold of dimension $2d$, endowed with the
effective Hamiltonian action of the torus $\R^d/L$ such that the image of the corresponding moment mapping 
is exactly $P$. Namely $M$ is the symplectic toric manifold $(M_P,\R^d/L,\omega,\Phi)$ corresponding to $P$. 
\end{prop}
Notice that there are no conditions on the vectors $v_j$, $j=h+1,\ldots, n$. They only have to lie in the lattice $L$. As in the previous sections, repetitions are allowed, and even zero vectors. These last yield degenerate half-spaces, coinciding with the whole space $(\R^d)^*$. 
We are allowed to add as many half spaces as we want. The proof of the Delzant procedure applies with no substantial changes. The key observation is that $\text{N}$ acts freely on $\Psi^{-1}(0)$. 
We will later need the explicit expression of $\Psi$: let $M\in M_{n\times (n-d)}(\R)$ be a matrix whose columns give a basis of $Rel(V)$. Then, expressed in components with respect the basis of $(Rel(V))^*$ dual to that basis, we have: 
\begin{equation}\label{momentmap}\Psi(\vz)=(|z_1|^2+l_1,\ldots,|z_n|^2+l_n)M.\end{equation}  
In particular, when the configuration is balanced, we can always take the first column to be $(1,\ldots,1)$, so that one of the components of $\Psi$ is 
\begin{equation}\label{sphere}(\sum_{i=1}^n(|z_i|^2+l),\end{equation} with $l=\sum_{i=1}^nl_i$. Therefore
$\Psi^{-1}(0)$ is contained in the sphere $S^{2n-1}_{r}$ of radius $r=-l$.
\begin{rem}\label{liscio}{
If, for some indices $\{j_1,\ldots,j_r\}\subset\{h+1,\ldots,n\}$, the third requirement is dropped, that is the hyperplanes
$\{\mu\in (\R^d)^*\;|\;\langle\mu,v_{j_k}\rangle=l_{j_k}\}$ {may} intersect $P$, the construction can still be followed step by step. But, in this case, $0$ is not a regular value of $\Psi$ and the level set $\Psi^{-1}(0)$ is not smooth. However, the reduced space continues to be smooth. This phenomenon was observed and thoroughly investigated by Guillemin-Sternberg in \cite{gs1}. Since the Delzant condition ensures that for each  $\vz\in\Psi^{-1}(0)$ the isotropy group $N_{\vz}$ is connected, there are no orbifold singularities in the quotient $\Psi^{-1}(0)/\text{N}$.}
\end{rem}
Recall that the Delzant procedure was generalized by Lerman-Tolman \cite{LT} to the cases in which $P$ is a  
simple convex polytope, rational with respect to a lattice $L$, and the vectors $v_i$ are not necessarily primitive in $L$. Prato further generalized the Delzant construction, so as to include the nonrational setting:
let $P$ be simple convex polytope. A triple $(P,Q,\{v_1,\ldots,v_h\})$, with $Q$ a quasilattice and
$v_i\in Q$ generators of the normal fan to $P$, will be called a symplectic Prato's datum. The resulting
reduced space is a symplectic toric quasifold $M$, of dimension $2d$, determined by the triple 
$(P,Q,\{v_1,\ldots,v_h\})$.
When $Q=L$ is a lattice, the Lerman-Tolman case is recovered, the resulting reduced space is the
symplectic toric orbifold determined by $(P,Q,\{v_1,\ldots,v_h\})$.  
When $Q=L$ is a lattice, the vectors are primitive and the polytope is Delzant, the resulting reduced space is the symplectic toric manifold $M_P$.

We now state a proposition similar to Proposition~\ref{sopra} for the above described cases. This variant
may be useful for all
cases such that
$$\text{Span}_\Z\{v_1,\ldots,v_h\}\subsetneqq Q,$$ where
$Q$ is a (quasi)lattice. 
This can happen when the polytope is not Delzant or when the first $h$ vectors are not primitive 
($Q=L$ is a lattice) or in the nonrational setting.
We then have $\exp(Rel(v_1,\ldots,v_h))\subsetneqq\text{N}$. That is $\text{N}$ is not connected.
 
By extending Proposition~\ref{sopra} we obtain the same reduced spaces resulting from the generalized Delzant procedure,
however, the group $\text{N}$ with respect to which we perform the symplectic reduction is connected.
It is enough to increase the number of half-spaces, exactly as observed in Section~\ref{Toric and LVMB}, in the complex set-up.
\begin{prop}\label{estensione} Let $(P,Q,\{v_1,\ldots,v_h\})$ a symplectic Prato's datum. Let $n\geq h$ and let $P$ be given by
$$P=\cap_{j=1}^{n}\{\mu\in(\R^d)^*\;|\;\langle\mu,v_j\rangle\geq l_j\},$$
where:
\begin{itemize}
\item $\text{Span}_\Z\{v_1,\ldots,v_n\}=Q$;
\end{itemize} 
and, if $n>h$
\begin{itemize}
\item $P\subset\{\mu\in (\R^d)^*\;|\;\langle\mu,v_j\rangle>l_j\},$ with $j=h+1,\ldots,n$.
\end{itemize}
Let $V=(v_1,\ldots,v_n)$ and consider the subgroup
$\text{N}=\exp(Rel(V))$ in $(S^1)^n$. Let $\Psi$ be the moment mapping with respect to
the induced action of $\text{N}$ on $\C^{n}$. Then the reduced
space $M=\Psi^{-1}(0)/\text{N}$ is endowed with the effective Hamiltonian action of the quasitorus $\R^d/Q$.
The image of the corresponding moment mapping $\Phi$ is exactly $P$. Moreover
$M$ is equivariantly symplectomorphic to
the sympletic quasifold (orbifold) determined by $(P,Q,\{v_1,\ldots,v_h\})$.
\end{prop}
The generalized Delzant procedure \cite{p} applies with no essential changes. We only have to check that 
$M$ can be identified with the symplectic quasifold $M'$ determined by $(P,Q,\{v_1,\ldots,v_h\})$.
We outline the argument: let $M'=(\Psi')^{-1}(0)/\text{N'}$ be the symplectic quasifold corresponding to
$(P,Q,\{v_1,\ldots,v_h\})$. Let $\Phi'$ be the moment mapping with respect to action of the quasitorus $\R^d/Q$ such that
$\Phi'(M')=P$. Consider the natural inclusion $\R^h\subset\R^n$.
We may view $Rel(v_1,\ldots,v_h)$ as a subset of $Rel(V)$. This allows to write explicitly an equivariant map 
$(\Psi')^{-1}(0)\stackrel{i}{\rightarrow}{\Psi}^{-1}(0)$ that induces a symplectomorphism $M'\rightarrow M$.
The key point is to verify that $\text{N}'$ at $\vz$ and $\text{N}$ at $i(\vz)$ have the same stabilizers, where
$\text{N}'$ and $\text{N}$ are defined in (\ref{enne}) and, in particular, $\text{N}=\exp(Rel(V))$.

{Finally notice that, as in Remark~\ref{liscio}, we can consider the cases in which the condition $P\subset\{\mu\in (\R^d)^*\;|\;\langle\mu,v_j\rangle>l_j\}$ is dropped for some indices in $\{h+1,\ldots,n\}$. Although the level set
is singular, the resulting reduced space is a manifold, orbifold or quasifold. 
These degenerate cases turn out to be relevant in various instances, see for example \cite{alf}, \cite[Remark~2.6]{quasicut} and \cite{quasireduction}.}
\begin{rem} {Consider the datum $(P,Q,\{v_1,\ldots,v_h\})$ and a non minimal, non degenerate presentation of $P$. 
By applying the generalized Delzant procedure to this presentation, as in Proposition~\ref{estensione}, we obtain the  reduced space corresponding to $(P,Q,\{v_1,\ldots,v_h\})$. However, the group that we use for the reduction is always connected and of arbitrarily high dimension. In the Delzant case, by Remark~\ref{liscio}, non degeneracy can be dropped.}
\end{rem}
Symplectic and complex quotients can often be identified via a natural homeomorphism or
diffeomorphim. These kind of results can be dated back to the work 
by Kempf and Ness \cite{kempfness} and Kirwan \cite{kirwan}, and apply to a number of
settings, of finite and infinite dimension.
A model example is given by the standard actions of $S^1$ and $(S^1)_\C=\C^*$ on $\C^{n}$. It is a toric example: let the polytope $P$ be the simplex in $\R^{n-1}$, with vertices the origin and $re_1,\ldots,re_{n-1}$, with $r=-l$. The fan $\Delta$ is its normal fan.
By (\ref{sphere}) the moment mapping
with respect to the $S^1$-action is $\Psi(\vz)=\sum_{i=1}^{n}|z_i|^2-r$. The zero level set is therefore the sphere
$S^{2n-1}_r$. While $U'(\Delta)=U(\bT)=\C^{n}\setminus\{0\}$.
We have the following diagram:
$$\xymatrix{
S^{2n-1}_r\ar@{^{(}->}[r]\ar[d]&\C^{n}\setminus\{0\}\ar[d]\\
S^{2n-1}_r/S^1\ar[r]^{\chi}&\C P^{n-1}
}$$
The inclusion of the zero level set in $U'(\Delta)$ induces a diffeomorphim at the level of the quotients, moreover the
symplectic structure on $S^{2n-1}_r/S^1$ is compatible with the complex structure of $\C P^{n-1}$ and induces a 
K\"ahler structure on the projective space, multiple of the Fubini-Study metric.
This principle holds also in our case.
Assume from now on the hypotheses of Proposition~\ref{estensione} and take
$\{V,\bT\}$ to be the corresponding triangulated vector configuration, namely $V=\{v_1,\ldots,v_n\}$ and $\bT$ the triangulation determined by the fan $\Delta$. A simple adaptation of \cite[Theorem~3.2]{cx} yields the following 
\begin{prop}\label{thechimap}
The level set $\Psi^{-1}(0)$ is contained in $\hat{U}(\bT^*)$. This inclusion induces an equivariant diffeomorphism
$$\chi\colon \Psi^{-1}(0)/\exp(Rel(V))\rightarrow\hat{U}(\bT^*)/\exp(Rel(V)_\C)$$
with respect to the actions of the quasitorus $\R^d/Q$ and the complex quasitorus $\C^d/Q$. Moreover
the induced symplectic form on the complex quasifold $\hat{U}(\bT^*)/\exp(Rel(V)_\C$ is K\"ahler. 
\end{prop}
We now insert our foliated complex manifold $(N,\calF)$ in this picture. Let $\{\Lambda,\bT^*\}$ an LVMB datum Gale dual to $\{V,\bT\}$. Then the columns of the matrix
$$
M=\begin{bmatrix}
1\horiz\ \Lambda^\R_1\ \horiz \\
      \vdots\ \  \\
1\horiz\ \Lambda^\R_{n}\ \horiz \\
\end{bmatrix}$$
are a basis of $Rel(V)$. Recall formulae (\ref{momentmap}) and (\ref{sphere})
and consider the following commutative diagram:
$$\xymatrix{
S_r^{2n-1}\supset \Psi^{-1}(0)\ar@{^{(}->}[r]^{\hat{i}}\ar[d]_{S^1}
\qquad&\qquad\hat{U}(\bT^*)\subset\C^n\setminus\{0\}\ar[d]^{\C^*}
\\
S^{2n-1}_r/S^1\supset \Psi^{-1}(0)/S^1\ar@{^{(}->}[r]^{i}\ar[dd]_{\exp({Rel(V)})/S^1}\ar[dr]^f\qquad&\qquad U(\bT^*)\subset \C\P^{n-1}\ar[d]^{\exp(\h)}\\
&N\ar[d]^{\exp(\overline{\h})}\\
M\ar[r]^{\chi}&X}.
$$
The composition of the vertical maps on the left-hand side give the quotient map by the group $\exp(Rel(V))$, while
the composition of the vertical maps on the right-hand side give the quotient map by the group $\exp(Rel(V)_\C$.

{From Proposition~\ref{thechimap} and the proof of \cite[Theorem~3.2]{cx} we know that, for each $\vz\in\Psi^{-1}(0)\subset \hat{U}(\bT^*)$,
the orbit $\exp(Rel(V)_{\C})\cdot\vz$ intersects $\Psi^{-1}(0)$ exactly in the orbit $\exp(Rel(V))\cdot\vz$. Moreover
the orbit  $\exp(iRel(V))\cdot\vz$ intersects $\Psi^{-1}(0)$ transversally, only in $\vz$. This still holds for $\Psi^{-1}(0)/S^1$ and the
orbits $\exp(Rel(V)_{\C}/\C)\cdot\vz$ and $\exp(iRel(V)/i\R)\cdot\vz$. Since $\exp(\h)\cap (S^1)^n=\{1\}$
we have that the orbit $\exp(\h)\cdot\vz$ intersects $\Psi^{-1}(0)/S^1$ transversally, only in $\vz$.
Therefore $f$ is a diffeomorphism. This is similar to \cite[Pag~83]{M}. Furthermore, $f$ sends
$\exp(Rel(V))/S^1$-orbits in $\Psi^{-1}(0)/S^1$
onto $\exp(\overline{h})$-orbits in $N$.}
Now let $\omega$ be the presymplectic form on $\Psi^{-1}(0)/S^1$
induced by the standard K\"ahler form in $\C^n$ (cf. \cite[Proposition~D]{MV}). Then, by the symplectic reduction properties, at each point $p\in \Psi^{-1}(0)/S^1$, the kernel of $\omega$ is the tangent space to the $\exp(Rel(V))/S^1$--orbit through $p$. Thus the kernel of $(f^{-1})^*(\omega)$ is the tangent space to the leaf $\exp(\overline{h})\cdot f(p)$. Moreover, $(f^{-1})^*(\omega)$  is clearly compatible with the complex structure on $N$. Therefore the form $(f^{-1})^*(\omega)$ endows $(N,\calF)$ with a transversely K\"ahler structure. This is of course, as we have seen in Section~\ref{polytopal}, a well known result, that we recover in our variant of the Delzant construction picture. This responds to a suggestion of \cite[pg.~74]{MV}.

We conclude by recalling some recent results of Ishida \cite{I2}, where he adopts symplectic methods to the study of 
complex manifolds admitting a maximal torus action. He proves a theorem analogous to the Atiyah and Guillemin-Sternberg convexity theorem, in the context of foliated manifolds endowed with a transversely symplectic form. 
Using this theorem he proves the direct implication of Theorem~\ref{transverse}, conjectured in \cite{CZ}.

Finally, in connection with this symplectic view-point, let us mention the very recent preprint by Nguyen and Ratiu \cite{NR}.
 
\medskip

{\em Acknowledgements} {The authors would like to thank Leonor Godinho for pointing out the relevance of \cite{gs1}}.

\bigskip

{\small 

\noindent Fiammetta Battaglia, Universit\`a di Firenze, Dipartimento di Matematica e Informatica U. Dini, Via S. Marta 3, 50139 Firenze, Italy.

\noindent
e-mail: fiammetta.battaglia@unifi.it}

\medskip

{\small 
\noindent
Dan Zaffran, College Of Marin,\\ 
835 College Ave, Kentfield, CA 94904

\noindent
e-mail: dan.zaffran@gmail.com}

\end{document}